\documentclass[12pt]{amsart}
\usepackage{amssymb}
\usepackage{graphicx}
\newtheorem{thm}{Theorem}[section]
\newtheorem{cor}[thm]{Corollary}
\newtheorem{lem}[thm]{Lemma}

\newtheorem{qu}[thm]{Question}
\theoremstyle{definition}

\numberwithin{equation}{section}

\begin{document}
\title[Engel graph]{Engel graph associated with a group}%
\author{Alireza Abdollahi }%
\address{Department of Mathematics, University of Isfahan, Isfahan 81746-73441, Iran}%
\email{a.abdollahi@math.ui.ac.ir}%
\thanks{This research  was in part supported by a grant from IPM (No. 86200021).}%
\begin{abstract}
Let $G$ be a non-Engel group and let $L(G)$ be the set of all left
Engel elements of $G$.
 Associate with $G$ a graph $\mathcal{E}_G$ as follows: Take $G\backslash L(G)$ as vertices of
$\mathcal{E}_G$ and join two distinct vertices $x$ and $y$
whenever $[x,_k y]\not=1$ and $[y,_k x]\not=1$ for all positive
integers $k$. We call $\mathcal{E}_G$, the Engel graph of $G$. In
this paper we study the graph theoretical properties of
$\mathcal{E}_G$.
\end{abstract} \maketitle
\section{\bf Introduction}
Let $G$ be a group and $x_1,\dots,x_n\in G$. For all $n>0$ we
define inductively $[x_1,\dots,x_n]$ as follows: $[x_1]=x_1$ and
$$[x_1,\dots,x_n]
=[x_1,\dots,x_{n-1}]^{-1}x_n^{-1}[x_1,\dots,x_{n-1}]x_n
\;\;\text{for all}\;\; n>1.$$ If $x_2=\dots=x_n$, then we denote
$[x_1,\dots,x_{n}]$ by $[x_1,_{n-1} x_2]$. Note that
$[x_1]=[x_1,_0 x_2]=x_1$ and $[x_1,x_2]=x_1^{-1}x_2^{-1}x_1x_2$.\\
An element $x$ of $G$ is called {\em left Engel} if for every
element $a \in G$, there exists a positive integer $k$ such that
$[a,_k x]=1$. If the integer $k$ is fixed for any element $a$,
then the element $x$  is called {\em left $k$-Engel}.  An element
$x$ is called {\em bounded left Engel} if it is {\em left
$k$-Engel} for some positive integer $k$. The sets of all left
Engel  and bounded left Engel elements of $G$ are denoted by
$L(G)$ and $\overline{L}(G)$, respectively. A group $G$
is called an {\em Engel} group, if $L(G)=G$.\\
 Associate with a non-Engel group $G$ a (simple) graph
$\mathcal{E}_G$  as follows: Take $G\backslash L(G)$ as vertices
of $\mathcal{E}_G$ and join two distinct vertices $x$ and $y$
whenever $[x,_k y]\not=1$ and $[y,_k x]\not=1$ for all positive
integers $k$. We call $\mathcal{E}_G$, the Engel graph of $G$.
\begin{center}
\includegraphics{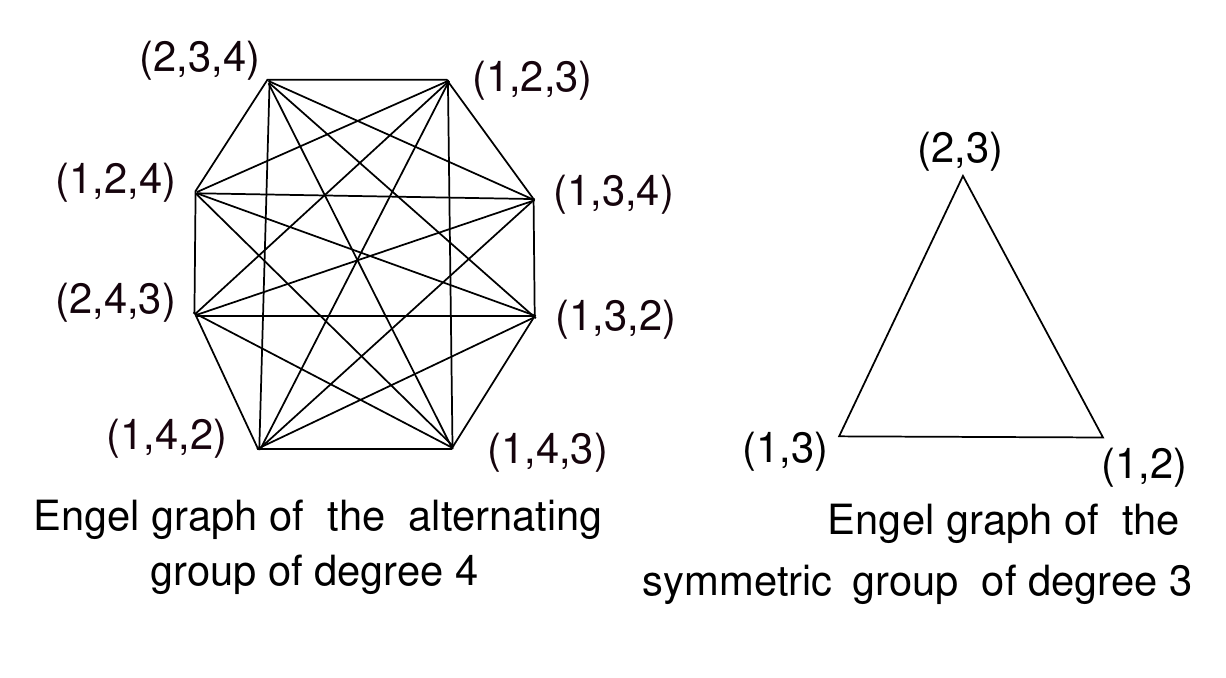}\\
Figure 1
\end{center}
In this paper we study the graph theoretical properties of
$\mathcal{E}_G$. One of our motivations for associating with a
group such kind of graph is a problem posed by Erd\"os: For a
group $G$, consider a graph $\Gamma$ whose vertex set is $G$ and
join two distinct elements if they do not commute. Then he asked:
is there a finite bound for the cardinalities of cliques  in
$\Gamma$, if $\Gamma$ has no infinite clique? (By a clique of a
graph $\Delta$ we mean a set of  vertices of $\Delta$ which are
pairwise adjacent. The largest size (if it exists) of cliques
 of $\Delta$ is called the clique number of $\Delta$; it will be denoted by $\omega(\Delta)$).\\
Neumann \cite{Neu} answered positively  Erd\"os' problem  by
proving that such groups are exactly the centre-by-finite groups
and the index of the centre can be considered as the requested
bound in the problem.\\
In fact,  groups with some conditions on  cliques of their Engel
graphs  have been already studied, without explicitly specifying
that such a graph had been under consideration.  For example,
Longobardi and Maj \cite{LonMaj} (also Endimioni \cite{End})
proved that if $G$ is a finitely generated soluble group in which
every infinite subset contains two distinct elements $x$ and $y$
such that $[x,_k y]=1$ for some integer $k=k(x,y)$, then $G$ is
finite-by-nilpotent. The following result is an easy consequence
of the latter  which may be considered as an answer to a Erd\"os
like question  on Engel graphs.
\begin{thm}
If $G$ is a non-Engel finitely generated soluble group, then
$\mathcal{E}_G$ has no infinite clique if and only if the clique
number of $\mathcal{E}_G$ is finite.
\end{thm}
In \cite{A1}, the condition $\mathcal{E}(n)$ ($n\in\mathbb{N}$) on
groups was considered by the author: Every subset of elements
consisting of more than $n$ elements possesses a pair $x,y$ such
that $[x,{}_ky]=1$ for some $k=k(x,y)\in\mathbb{N}$. This simply
means  that we were studying groups $G$ such that every clique of
$\mathcal{E}_G$ has size at most $n$. It is shown in \cite{A1}
that all finite groups and all finitely generated soluble groups
satisfying $\mathcal{E}(n)$ are nilpotent for $n\leq 2$ and that
all finite groups in $\mathcal{E}(n)$ are soluble for $n\leq 15$.
Therefore we can summarize the latter results  in terms  of Engel
graph as following.
\begin{thm}\label{thm1}
Let $G$ be a finite or a finitely generated soluble non-Engel
group. Then $\omega(\mathcal{E}_G)\geq 3$. If $G$ is finite and
$\omega(\mathcal{E}_G)\leq 15$, then $G$ is soluble.
\end{thm}
In \cite{A2}, it is proved that if $G$ is a finitely generated
soluble group satisfying $\mathcal{E}(n)$, then the index of the
hypercentre of $G$ is bounded by a function of $n$. Also it is
proved that if $G$ is a finite  group satisfying $\mathcal{E}(n)$,
then the index of the Fitting subgroup of $G$ is bounded by a
function of $n$. Note that if $G$ is  either a finite group or a
soluble group, then $L(G)$ is a subgroup and coincides with the
Fitting subgroup of $G$ (see \cite{Baer} and \cite[Proposition
3]{Gr}).
\begin{thm}
Let $G$ be a finite or a finitely generated soluble non-Engel
group such that $\omega(\mathcal{E}_G)$ is finite. Then the index
of the Fitting subgroup of $G$ is finite and bounded by a
function  of $\omega(\mathcal{E}_G)$.
\end{thm}
Also in \cite{A2} finite centerless groups $G$ satisfying the
condition $\mathcal{E}(n)$ for $n\leq 15$ are characterized. This
of course implies that  we have  a characterization of all finite
centerless groups $G$ with $\omega(\mathcal{E}_G)\leq 15$.\\
In section 2, we study the connectedness of $\mathcal{E}_G$ for a
non-Engel group $G$. We do not know whether there is a non-Engel
group $G$ such that  $\mathcal{E}_G$ is not connected. We give
some classes of groups whose Engel graphs are connected.
Throughout section 2 we generalize some known results on the set
of left Engel elements by defining  new types of Engel elements.
Most of the results are about the connectedness  of certain
subgraphs of a non-Engel group.  In section 3, we characterize
finite groups whose Engel graphs are planar. In section 4 we
shortly study  groups with isomorphic Engel graphs and we show
that the Engel graph of a finite group cannot be isomorphic to
the Engel graph of an infinite one.
\section{\bf Connectedness of Engel graph }
In this section, we study the connectedness of the Engel graph
$\mathcal{E}_G$ for a non-Engel group $G$. Before starting to show
the results, we recall some concepts for a simple graph $\Delta$.
A {\it path} $P$ in $\Delta$  is a sequence $v_0-v_1-\cdots-v_k$
whose terms are vertices of $\Delta$ such that for any
$i\in\{1,\dots,k\}$, $v_{i-1}$ and $v_i$ are adjacent.  In this
case $P$ is called a path between $v_0$ and $v_k$. The number $k$
is called the {\it length} of $P$.
 If $v$ and $w$ are vertices in $\Delta$, then by definition $d(v,v)=0$ and
 whenever  $v\not=w$,
$d(v,w)$ denotes the length of  the shortest path between $v$ and
$w$ if a path exists, otherwise $d(v,w)=\infty$ and we call
$d(v,w)$ the distance between $v$ and $w$ in $\Delta$. We say that
$\Delta$ is {\it connected} if there is a path between each pair
of  vertices of $\Delta$. If $\Delta$ is connected, then the
largest distance between all pairs of the vertices of $\Delta$ is
called the {\it diameter} of $\Delta$, and
it is denoted by $\text{diam}(\Delta)$.\\

Some results and examples suggest that the Engel graph of any
finite non-Engel group must be connected. We have checked the
Engel graphs of all non-Engel finite groups of orders at most 570
(except ones of order 384!) by GAP \cite{GAP} using the package
GRAPE and the following programme written in GAP. All these
groups have connected Engel graphs  with diameter 1 or 2.
\begin{verbatim}
c:=function(a,b,G) local n,L,s,i,A,M; n:=Size(G); M:=1; s:=0; if
Comm(a,b)=Identity(G) then s:=1; fi; A:=a;
 while s=0 do
A:=Comm(A,b); M:=M+1; if A=Identity(G) then s:=1; fi; if M>n and
A<>Identity(G) then s:=-1; fi; od; if s=0 or s=-1 then L:=true;
fi; if s=1  then L:=false; fi; return L; end;
c2:=function(a,b,G)
return c(a,b,G) and c(b,a,G); end;
LoadPackage("GRAPE");
EngelGraph:=function(G)
return
Graph(G,Difference(G,FittingSubgroup(G)),OnPoints,function(x,y)
return c2(x,y,G); end); end;
DiametersOfEngelGraphs:=function(n,m) return
Union(List([n..m],j->List(AllSmallGroups(j,IsNilpotent,false),
i->Diameter(EngelGraph(i))))); end;
\end{verbatim}
The difficulty to settle the question of whether Engel graphs of
finite non-Engel groups are connected, is that the relation $\sim$
on $G$ defined as
$$x \sim y \Leftrightarrow [x,_k y]=1 \;\text{for some}\; k\in\mathbb{N}$$
is not symmetric, for example in the symmetric group of degree 3,
we have  $(1,2)\sim(1,2,3)$ but $(1,2,3)\not\sim(1,2)$. Indeed, if
it is true that $\mathcal{E}_G$ is connected, it should be true
that $\mathcal{E}_G$ has no isolated vertex, i.e. a vertex $v$
such that $d(v,w)=\infty$ for every vertex $w\not=v$. We were
unable to prove the latter for all groups, but we shall prove it
for some wide classes of groups. Our results generalize some
well-known results on left Engel elements and Engel sets, as we
consider a new type of Engel elements in groups. In the
definition of edges of the complement of an Engel graph, since we
want to have a simple graph, we use the symmetrized of the
relation $\sim$, i.e. $$x \sim' y \Leftrightarrow \text{either}\;
[x,_k y]=1 \;\text{or}\; [y,_k x]=1 \;\text{for some}\;
k\in\mathbb{N}.$$ This motivates to consider the following types
of Engel conditions on elements of a
group.\\
Let $G$ be a group.  An element $a$ of $G$ is called {\em randomly
power Engel} if for every $g\in G$, there exists a sequence
$a_1,\dots,a_k$ of elements of $\left<a\right>$ with
$\left<a\right>=\left<a_i\right>$ for all $i\in\{1,\dots,k\}$ such
that  either $[a_1^g,a_2,\dots,a_k]=1$ or
$[a_1,a_2^g,\dots,a_k^g]=1$. The word ``randomly'' has been
selected for the name of such an element, because of the
``either...or'' sentence involved in the definition and the word
``power'' is for the fact that, the elements $a_1,\dots,a_n$ are
indeed powers of $a$. If the integer $k$ in above is the same for
all $g\in G$, then we say that $a$ is a {\em randomly power
$k$-Engel} element. If all the elements $a_1,\dots,a_k$ can be
always chosen equal to $a$, then $a$ is called a {\em randomly
Engel element} of $G$, {\em randomly $k$-Engel} if $k$ is the
same for all $g\in G$. An element $a$ of $G$ is called {\em
bounded randomly {\rm (}power{\rm)} Engel} if it is randomly
(power) $k$-Engel for some positive integer $k$. The set of
randomly power Engel elements in a group $G$ is denoted by
$\mathcal{L}(G)$ and we denote by $\overline{\mathcal{L}}(G)$ the
set of bounded randomly power Engel elements of $G$.  Clearly we
have $L(G)\subseteq \mathcal{L}(G)$, $\overline{L}(G)\subseteq
\overline{\mathcal{L}}(G)$ and these  sets are invariant under
conjugation of elements of $G$.\\
We need the following lemma in the sequel.
\begin{lem}\label{lem0}
Let $a$ and $g$ be elements of a group such that the normal
closure of $a$ in $\left<a,g\right>$ is  abelian. If
$[a,g^{t_1},\dots,g^{t_k}]=1$ for some $t_1,\dots,t_k
\in\mathbb{N}$, then $[a,_k g^{m}]=1$, where $m$ is any positive
integer divisible by the least common multiple of  $t_1,\dots,
t_k$.
\end{lem}
\begin{proof}
Since $\left<a\right>^{\left<a,g\right>}$ is abelian, we may write
\begin{align*}
[a,g^{m},g^{t_2},\dots,g^{t_k}]=[[a,g^{t_1}]^{(g^{t_1})^{\frac{m}{t_1}-1}+(g^{t_1})^{\frac{m}{t_1}-2}+\cdots+g^{t_1}+1},g^{t_2},\dots,g^{t_k}]&\\
[a,g^{t_1},g^{t_2},\dots,g^{t_k}]^{(g^{t_1})^{\frac{m}{t_1}-1}+(g^{t_1})^{\frac{m}{t_1}-2}+\cdots+g^{t_1}+1}=1.&
\end{align*}
Now since the normal closure of $[a,g^{m}]$ in $\left<a,g\right>$
is also abelian, by induction on $k$, we have that $[a,_k g^m]=1$.
This completes the proof.
\end{proof}
\begin{lem}\label{lem1}
Let $A$ be a normal abelian   subgroup of a group $G$ and let
$g\in G$.
\begin{enumerate}
\item If  $g\in \mathcal{L}(A\left<g\right>)$, then
$A\left<g\right>$ is locally nilpotent. \item If $g$ is a randomly
power $k$-Engel element of $A\left<g\right>$, then
$A\left<g\right>$ is nilpotent of class at most $k+1$.
\end{enumerate}
\end{lem}
\begin{proof}
Let $K=A\left<g\right>$. To prove (1) and (2) it is enough to show
that $g\in L(K)$ and
$g$ is a left $k$-Engel element, respectively.\\
  Let $a\in A$, $k\in\mathbb{N}$ and
$t_1,t_2,\dots,t_k \in \mathbb{Z}$. Then, since $A$ is a normal
abelian subgroup of $G$, we can write
\begin{align*}
[(g^{t_1})^a,g^{t_2},\dots,g^{t_k}]=[[g^{t_1},a],g^{t_2},\dots,g^{t_k}]=& \\
[[a,g^{t_1}]^{-1},g^{t_2},\dots,g^{t_k}]=[a,g^{t_1},g^{t_2},\dots,&g^{t_k}]
^{-1},
\end{align*} and
\begin{align*}
[g^{t_1},(g^{t_2})^a,\dots,(g^{t_k})^a]=[(g^{t_1})^{a^{-1}},g^{t_2},\dots,g^{t_k}]^{a}=& \\
[[g^{t_1},a^{-1}],g^{t_2},\dots,g^{t_k}]^a=[a,g^{t_1},g^{t_2},\dots,g^{t_k}].
\end{align*}
Since $g\in \mathcal{L}(K)$, it follows from these equalities that
for any $a\in A$, there exists a sequence $t_1,\dots,t_k$ of
integers with $\left<g^{t_i}\right>=\left<g\right>$ for every
$i\in\{1,\dots,k\}$ such that
$$[a,g^{t_1},g^{t_2},\dots,g^{t_k}]=1. \eqno{(*)}$$ Now if $g$ is of infinite order,  then $t_1,\dots,t_k
\in\{1,-1\}$ and by $(*)$ it is easy to see that  for any $a\in
A$, there exists a positive integer $k$ such that $[a,_k g]= 1$.
Since $A$ is abelian and normal in $K$, it follows that $g \in
L(K)$. Now assume that $g$ is of finite order and let $m$ be the
product of all positive integers $t\leq |g|$ such that
$\gcd(t,|g|)=1$. It follows from $(*)$ and Lemma \ref{lem0} that
$[a,_k g^m]=1,$ which implies  that $g^m \in L(K)$. Now by  a
result of Gruenberg \cite[Proposition 3]{Gr}, $L(K)$ is a
subgroup of $G$ and so $g\in L(K)$. This shows that, in any case
$g\in L(K)$. The above argument also shows that if $g$ is a
randomly power $k$-Engel element of $K$, then $g$ is a left
$k$-Engel element of $K$. This completes the proofs of (1) and
(2).
\end{proof}
\begin{thm}\label{Nil-Ab}
 Let $G$ be a nilpotent-by-abelian non-Engel group. If $\mathcal{E}_{\frac{G}{H'}}$ is connected for some nilpotent
subgroup $H$ containing $G'$, then  $\mathcal{E}_G$ is connected.
 Moreover, in this case, the diameter
of $\mathcal{E}_G$ is at most
$\max\{\text{diam}(\mathcal{E}_{\frac{G}{H'}}),2\}$.
\end{thm}
\begin{proof}
For every $a\in G$, $a\in L(G)$ if and only if $aH'\in
L(\left<a\right>H/H')$; for if $aH'\in L(\left<a\right>H/H')$,
then $\left<a\right>H/H'$ is locally nilpotent by Lemma
\ref{lem1}. Since $H$ is nilpotent, a Hall-type result of Plotkin
\cite{plot} (see also \cite{Rob4}) implies that $\left<a\right>H$
is locally nilpotent. Since $G'\leq H$, $\left<a\right>^G\leq
\left<a\right>H$ and so $a\in L(G)$.
 It follows that $xH'$ is a vertex of
$\mathcal{E}_{\frac{G}{H'}}$ if and only if $x$ is a vertex of
$\mathcal{E}_{G}$. Now let $x$ and $y$ be two distinct vertices in
$\mathcal{E}_G$. If $xH'=yH'$, then  there exists $z\in G$ such
that $zH'$ is an adjacent vertex to $xH'$. It follows that
$x-z-y$ is a path of length 2 between $x$ and $y$ in
$\mathcal{E}_G$. Now if $xH'\not=yH'$, then there exists a path
$P$ of length $d\leq \text{diam}(\mathcal{E}_{\frac{G}{H'}})$
between $xH'$ and $yH'$ in $\mathcal{E}_{\frac{G}{H'}}$. Now
 any set of preimages of the vertices of $P$ in $G$ (under the natural
epimorphism $G\rightarrow \frac{G}{H'}$)  forms  a path of length
$d$ between $x$ and $y$ in $G$. This completes the proof.
\end{proof}
\begin{thm}\label{nil-cyc}
If $G$ is a nilpotent-by-cyclic non-Engel group, then
$\mathcal{E}_G$ is connected and its  diameter is at most $6$. If
$G$ is nilpotent-by-abelian and $x\not\in L(G)$, then
 $\left<x\right>^G$ is non-Engel,
$\mathcal{E}_{\left<x\right>^G}$ is connected with  diameter  at
most $6$, and the induced subgraph of $\mathcal{E}_{G}$ on the
conjugacy class of $x$ in $G$ is connected with diameter at most
$2$.
\end{thm}
\begin{proof}
 Let $G=A\left<g\right>$, where $A$ is a normal nilpotent subgroup of $G$ and $g\in G$.
 By Theorem \ref{Nil-Ab} we may assume that  $A$
is abelian.   Let $g_1$ and $g_2$ be two
 distinct vertices of $\mathcal{E}_G$ such that $g_1=a_1g^{n}$
and $g_2=a_2g^{n}$ for some $a_1,a_2\in A$ and $n\in \mathbb{Z}$.
Since $A$ is a normal abelian subgroup of $G$,  $g_1$ is adjacent
to $g_2$ if and only if
 $$[a_1a_2^{-1},_k g^n]\not=1  \;\;\text{for all}\;\; k\in \mathbb{N}.$$ Since $g_1$ is a vertex of
$\mathcal{E}_G$, $g_1\not\in \mathcal{L}(G)$ by Lemma \ref{lem1}.
Thus $g_1$ is adjacent to  $g_1^y$ for some $y\in G$. Note that
$g_1^y=a_3 g^n$ for some $a_3\in A$. Therefore $[a_1a_3^{-1},_k
g^n]\not=1$ for all $k\in\mathbb{N}$. Now if $g_1$ and $g_2$ are
not adjacent, then $[a_1a_2^{-1},_m g^n]=1$ for some
$m\in\mathbb{N}$. Thus $[a_2a_3^{-1},_k g^n]\not=1$ for all $k\geq
m$ and so clearly we have that $[a_2a_3^{-1},_k g^n]\not=1$ for
all $k\in\mathbb{N}$. Hence $g_2$ is adjacent to $g_1^y$ and
$g_1-g_1^y-g_2$ is a path of length 2 between $g_1$ and $g_2$.
This implies that any two distinct vertices of $\mathcal{E}_G$ of
forms $a_1g^n$ and $a_2g^n$ are either adjacent or  there is a
path of length 2 between them so that the middle vertex in this
path can be   a suitable conjugate of either $g_1$ or $g_2$ (*).\\

Now we prove that if $a_1g^n$ is a vertex of $\mathcal{E}_G$ (with
$a_1\in A$ and $n\in\mathbb{Z}$), then there is a path of length
at most 3 between $a_1g^n$ and $g$. As $g^n\not\in\mathcal{L}(G)$
by Lemma \ref{lem1}, there exists $a\in A$ such that $(g^n)^a$ is
adjacent to $g^n$. This is equivalent to $[a,_k g^n]\not=1$ for
all $k\in\mathbb{N}$. Now Lemma \ref{lem0} implies that $[a,g^n,_k
g]\not=1$ and $[a,g,_k g^n]\not=1$ for all $k\in\mathbb{N}$.
Therefore $[(g^n)^a,_k g]\not=1$ and $[g,_k (g^n)^a]\not=1$ for
all $k\in\mathbb{N}$ and so $(g^n)^a$ is adjacent to $g$. By the
previous part, there is a path of length at most 2 between
$a_1g^n$ and $(g^n)^a$ and so there is a path of length 3 between
$a_1g^n$ and $g$.\\

Now let $a_2g^m$ be another vertex in $\mathcal{E}_G$, where
$a_2\in A$. Therefore by the latter paragraph, there is a path of
length $3$ between $a_2g^m$ and $g$. Thus $d(a_1g^n,a_2g^m)\leq
6$.  This completes
the proof of the first statement of the theorem.\\
Suppose now that $G$ is nilpotent-by-abelian and $x\not\in L(G)$.
Then $x\not\in L(\left<x\right>^G)$,  otherwise, by a result of
Gruenberg \cite[Proposition 3]{Gr}, $x$ lies in the Hirsch-Plotkin
radical of $\left<x\right>^G$, and since the latter subgroup is
normal in $G$, $x$ lies in the Hirsch-Plotkin radical of $G$,
which contradicts $x\not\in L(G)$. This implies that
$\left<x\right>^G$ is non-Engel. As $\left<x\right>^G\leq
\left<x\right>G'$, we have that $\left<x\right>^G$ is
nilpotent-by-cyclic, and so by the previous part
$\mathcal{E}_{\left<x\right>^G}$ is connected and
$\text{diam}(\mathcal{E}_{\left<x\right>^G})\leq 6$.\\
The last part of the theorem  easily follows  from (*).
\end{proof}
As far as we know, it is an open problem whether the set of left
Engel elements of an arbitrary group $G$ forms a subgroup. However
there are classes of groups $G$ in which not only $L(G)$ and
$\mathcal{L}(G)$ are  subgroups but also they are very
well-behaved. We shall see some of these classes in the following.
\begin{thm}\label{GR}
Let $G$ be a soluble group.
\begin{enumerate}
\item $\mathcal{L}(G)$ coincides with the Hirsch-Plotkin radical
of $G$ and is a Gruenberg group. In particular,
$\mathcal{L}(G)=L(G)$. \item $\overline{\mathcal{L}}(G)$ coincides
with the Baer radical of  $\; G$. In particular,
$\overline{\mathcal{L}}(G)=\overline{L}(G)$.
\end{enumerate}
\end{thm}
\begin{proof}
Let $g_1\in \mathcal{L}(G)$ and $g_2\in
\overline{\mathcal{L}}(G)$. By a result of Gruenberg
\cite[Proposition 3]{Gr} (see also \cite[12.3.3]{Rob}), it is
enough to show that $g_1\in L(G)$ and $g_2\in \overline{L}(G)$.
We argue by induction on the derived length $d$ of $G$. If $d\leq
1$, then $G$ is abelian and $g_1,g_2\in \overline{L}(G)$; thus we
can assume $d>1$ and write $A=G^{(d-1)}$. Now obviously $g_1A \in
L(G/A)$ and $g_2A \in \overline{L}(G/A)$. It remains to prove
that $g_1 \in L(\left<A,g_1\right>)$ and $g_2 \in
\overline{L}(\left<A,g_2\right>)$: this immediately follows from
Lemma \ref{lem1}.
\end{proof}
In \cite{Baer} Baer proved that in every group $G$ satisfying the
maximal condition on all subgroups, $L(G)$ coincides with the
Fitting subgroup of $G$, and  in \cite{Peng}, Peng  by using an
argument of Baer, generalized  the latter result as following.
\begin{thm} \label{Peng}
If $G$ is a group satisfying the maximal condition on abelian
subgroups, then $L(G)$   coincides with the Fitting subgroup of
$G$.
\end{thm}
Recently a weaker version of Theorem \ref{Peng} has been proved in
\cite[Theorem 1]{Mamontov}. In the following we generalize Theorem
\ref{Peng} to randomly power Engel elements.
\begin{thm}\label{main}
If $G$ is a group satisfying the maximal condition on abelian
subgroups, then $L(G)=\mathcal{L}(G)$. In particular,
$\mathcal{L}(G)$ coincides with the Fitting subgroup of $G$.
\end{thm}
\begin{proof}
Let $a\in \mathcal{L}(G)$. It is enough to show that
$\left<a\right>^G$ is nilpotent.\\
The proof is similar (mostly {\em mutatis mutandis}) to that given
for \cite[Part 2, Lemma 7.22]{Rob1}, but it needs some little
change. With the notation of that proof, one has to consider the
subgroups
$$W_U = \langle N_U(I)\cap a^G\rangle \;\;\text{and}\;\;
W_V=\langle N_V(I)\cap a^G\rangle,$$ where $a^G$ denotes the set
of all conjugates of $a$ in $G$, and fix elements $v\in
(N_V(I)\cap a^G)\backslash I$ and $u\in (N_U(I) \cap
a^G)\backslash I$. Then make use of Theorem \ref{GR} and of the
fact that if $\langle a^t\rangle= \langle a\rangle$ for some
integer $t$, then $\langle u^t\rangle=\langle u\rangle$ and
$\langle v^t\rangle=\langle v\rangle$.
\end{proof}
In \cite[Part 2, page 55]{Rob1}, a subset $S$ of a group $G$ is
called an {\em  Engel set} if given $x$ and $y$ in $S$ there is an
integer $n=n(x,y)$ such that $[x,_n y]=1$. As a corollary to
\cite[Part 2, Lemma 7.22]{Rob1} (Theorem \ref{Peng}, here), normal
Engel sets in a group $G$ satisfying the maximal condition on
abelian subgroups are characterized in \cite[Part 2, Theorem
7.23]{Rob1} as the subsets of the Fitting subgroup of $G$. Recall
that  a normal set in a group  is a set closed under conjugation.
We call a subset $R$ of a group $G$ a {\em randomly Engel set} if
given $x$ and $y$ in $R$ there is an integer $n=n(x,y)$ such that
either $[x,_n y]=1$ or $[y,_n x]=1$. Generally a randomly Engel
set is not an Engel set, consider for example $\{(1,2),(1,2,3)\}$
in the symmetric group of degree 3, but we shall see in the
following result  that the normal ones are the same in certain
groups.
\begin{thm}\label{maincor}
Let $G$ be a group satisfying the maximal condition on abelian
subgroups {\rm (}so  $G$ may be a finite group{\rm )}. Then a
normal subset of $G$ is a  randomly Engel set if and only if it is
contained in the Fitting subgroup of $G$. In particular,  an
element $x$ of $G$ lies in the Fitting subgroup of $G$ if and only
if for every $g\in G$, there exists a positive integer $k$  such
that either $[x^g,_k x]=1$ or $[x,_k x^g]=1$.
\end{thm}
\begin{proof}
Suppose that the normal subset $S$ of $G$ is a  randomly Engel set
and let $a\in S$. If $g\in G$, then $a^g\in S$ and so either
$[a^g,_n a]=1$ or $[a,_n a^g]=1$ for some integer $n$. It follows
that $a$ is a randomly Engel element of $G$ and so Theorem
\ref{main} implies that $\left<a\right>^G$ is nilpotent. Therefore
$S$ is contained in the Fitting subgroup of $G$. The converse is
clear.
\end{proof}
Now using Theorem \ref{maincor},  it is easy to  generalize
\cite[Satz 1]{Held} (see also \cite[Part 2, Theorem 7.24]{Rob1}).
The proof is mostly {\em mutatis mutandis}  the proof of
\cite[Part 2, Theorem 7.24]{Rob1} so we will omit it.
\begin{thm}\label{min}
Let $G$ be a group satisfying the minimal condition on subgroups
and suppose that the elements whose orders are powers of $p$ form
a randomly Engel set for each prime $p$. Then $G$ is a
hypercentral $\check{\text{C}}$ernikov group.
\end{thm}
\begin{cor}\label{min-max-sol}
Let $G$ be a  group in which  every two-generated subgroup of $G$
is either soluble or  satisfies the maximal condition on its
abelian subgroups. Then $L(G)=\mathcal{L}(G)$. In particular, if
$G$ is non-Engel, then $\mathcal{E}_G$ has no isolated vertex.
\end{cor}
\begin{proof}
The first part follows from Theorems \ref{GR} and \ref{main}.\\
If $G$ is non-Engel and $a$ is a vertex of $\mathcal{E}_G$, then
$a$ is not a randomly power  element of $G$. Thus there is a
conjugate of $a$ which is adjacent to $a$. This completes the
proof.
\end{proof}
The following result, which is of independent interest,
generalizes \cite[Theorem 1 (1)]{Ho}.
\begin{thm}
Let $G$ be a finite group and $p$ a prime number. Let $x$ be a
$p$-element  and $P$ a Sylow $p$-subgroup of $G$ such that for
every $y\in G$ with $x^y\in P$, the set $\{x,x^y\}$ is a randomly
Engel set. If $G$ is $p$-soluble, then $x\in P$.
\end{thm}
\begin{proof}
Suppose that $G$ is a counterexample of minimum order. Hence
$O_p(G)=1$ and $G=\left<P,x\right>$. Let $S$ be any minimal normal
subgroup of $G$. By minimality of $G$, we have that $G=PS$. Since
$G$ is $p$-soluble and $O_p(G)=1$, $S$ is a normal subgroup of
order relatively  prime to $p$. There exists $s\in S$ such that
$x^s\in P$. Thus $[x^s,_k x]=1$ or $[x,_k x^s]=1$ for some $k\in
\mathbb{N}$. Suppose that $[x^s,_k x]=1$ and $k$ is the least
non-negative integer with this property; obviously $k>0$. Then
$$[[x,s],_k x]=1 \eqno{(1)}$$ and $[[x,s],_{k-1} x]\not=1$.  Now we prove that if $k\geq 1$, then we lead  to a
contradiction.  Put $|x|=q$, then using $(1)$  we can write
$$1=[[x,s],_{k-1} x^q]=[[x,s],_{k-1}
x]^{q(x^{q-1}+\cdots+x+1)^{k-2}}= [[x,s],_{k-1}x]^{q^{k-1}}.$$
 Since
$[[x,s],_{k-1}x]\in S$ and $S$ is a $p'$-group, it follows that
$[[x,s],_{k-1}x]=1$ which  contradicts the minimality of $k$. This
completes the proof, in this case. If $[x,_k x^s]=1$, then
$[x^{s^{-1}},_kx]=1$ and   a similar argument completes the proof.
\end{proof}
\begin{thm}\label{dim=1}
Let $G$ be a non-Engel group.
\begin{enumerate}
\item If $L(G)$ is a subgroup of $G$ and $x$ is a vertex of
$\mathcal{E}_G$ such that $x$ is adjacent to every vertex of
$\mathcal{E}_G$, then $x^2=1$ and $C_G(x)=\left<x\right>$. \item
If  $L(G)$ is a subgroup of $G$, then {\rm
$\text{diam}(\mathcal{E}_G)=1$} if and only if $L(G)$ is a normal
abelian subgroup of $G$ without elements of order $2$ and for any
vertex $x$ of $\mathcal{E}_G$  we have $G=L(G)\left<x\right>$,
$L(G)\cap \left<x\right>=1$, $x^2=1$ and $g^x=g^{-1}$ for all
$g\in L(G)$.
 \item
If $L(G)$ is a subgroup of $G$ and  $G$ is periodic, then
$\mathcal{E}_G$ contains a vertex  which is adjacent to all other
vertices in $\mathcal{E}_G$ if and only if {\rm
$\text{diam}(\mathcal{E}_G)=1$.}
\end{enumerate}
\end{thm}
\begin{proof}
(1) \; Since $L(G)\leq G$, $x^{-1}$  is also a vertex and since
$x$ is adjacent to all vertices,  then $x=x^{-1}$. Now let $y\in
C_G(x)$ and $y\not=x$. Since $[x,y]=1$ and $x$ is adjacent to
every vertex, $y$ is not a vertex and so $y\in L(G)$. If  $y$ is
non-trivial, then $yx\not=x$ and by a similar argument, $yx\in
L(G)$. As $L(G)$ is a subgroup of $G$, it follows that $x\in
L(G)$, a contradiction. Hence $C_G(x)=\left<x\right>$.\\
(2) \; Let $\text{diam}(\mathcal{E}_G)=1$. By part (1),  $a^2=1$
for every vertex $a$ of $\mathcal{E}_G$. Now suppose $a$ and $b$
are two distinct vertices of $\mathcal{E}_G$. Then
$\left<ab\right>$ is a normal subgroup in the dihedral group
$\left<a,b\right>$ and so if $ab$ is  a vertex of
$\mathcal{E}_G$, then it is not adjacent to $a$. It follows that
$G/L(G)$ is of order 2 which implies that $G=L(G)\left<x\right>$
and $L(G) \cap \left<x \right>=1$, for every vertex  $x$  of
$\mathcal{E}_G$. If $g\in L(G)$, then $gx$ is a vertex of
$\mathcal{E}_G$. It follows that $(gx)^2=1$, or $g^x=g^{-1}$.\\
The converse in straightforward.\\
(3) \;   Let $x$ be a vertex of $G$ adjacent to all other vertices
of $\mathcal{E}_G$. By part (1), we have $C_G(x)=\left<x\right>$
and $x^2=1$. It follows that $\left<x\right>\cap
\left<x\right>^g=1$ for all $g\in G\backslash \left<x\right>$. Now
by Theorem 5 of \cite{Sh}, $A=G\backslash \{x^g \;|\; g\in G\}$ is
a normal abelian subgroup of $G$, such that $G=A\left<x\right>$
and obviously $A\cap \left<x\right>=1$. Thus $G$ is a solvable
periodic group which implies that  $G$ is locally finite. Let
$a\in A$ be a non-trivial element of $A$, then
$B=\left<a,a^x\right>$ is a finite abelian normal subgroup of $G$.
Hence $x$ induces    a fixed-point-free automorphism of order 2 in
$B$, which implies that $B$ is an abelian group of odd order  and
 $b^x=b^{-1}$ for all $b\in B$. Therefore $a^x=a^{-1}$ for all $a\in A$. Obviously we have that $L(G)=A$. Now part (2) completes
the proof.
\end{proof}
\begin{qu}
Is the hypothesis ``$L(G)\leq G$''  necessary in Theorem {\rm
\ref{dim=1}}?
\end{qu}
\begin{qu}
Let $G$ be a finite non-nilpotent group. Is it true that
$\mathcal{E}_G$ is connected? If so, is it true that {\rm
$\text{diam}(\mathcal{E}_G)\leq 2$}?
\end{qu}
\begin{qu}
In which classes of groups, the Engel graph of every non-Engel
group has no isolated vertex?
 \end{qu}
 \section{\bf Groups whose Engel graphs are planar}
A {\em planar graph} is a graph which can be drawn in the plane so
that its edges intersect only at end vertices.\\
 In this section we prove
 \begin{thm}
Let $G$ be a finite non-Engel group. Then $\mathcal{E}_G$ is
planar if and only if $G\cong S_3$, $D_{12}$ or $T=\left<x,y\;|\;
x^6=x^3y^{-2}=x^yx=1\right>$.
 \end{thm}
\begin{proof}
Suppose that $\mathcal{E}_G$ is planar. Then   $\mathcal{E}_G$ has
no subgraph isomorphic to $K_5$ (the complete graph with 5
vertices) or $K_{3,3}$ (the complete bipartite graph whose parts
have the same size 3) (see \cite[Corollary 4.2.11]{Di}). This
implies that $\omega(\mathcal{E}_G)\leq 4$. Now it follows from
Proposition 1.4 of \cite{A2} and Theorem \ref{thm1}, that
$\overline{G}=\frac{G}{Z^*(G)}\cong S_3$ or $A_4$, where $Z^*(G)$
is the hypercentre of $G$. As we see in Figure 1,
$\mathcal{E}_{A_4}$ has a subgraph isomorphic to $K_{3,3}$. Thus
$\overline{G}\cong S_3$. Now put $\bar{x}=xZ^*(G)$ for every $x\in
G$ and let $a,b\in G$ such that
$\overline{G}=\left<\bar{a},\bar{b}\right>\cong S_3$ where
$\bar{a}^3=\bar{b}^2=\bar{1}$. Then since every element of
$Z^*(G)$ is a right Engel element of $G$, we have that
$L(G)=\left<a\right>Z^*(G)$ and $$G\backslash L(G)=bZ^*(G) \cup
abZ^*(G) \cup a^2bZ^*(G).$$ Now because
$\{\bar{b},\bar{a}\bar{b},\bar{a}^2\bar{b}\}$ is a clique in
$\mathcal{E}_{\overline{G}}$ (see Figure 1) and every element of
$Z^*(G)$ is right Engel, we have that every element in $a^ib
Z^*(G)$ is adjacent to every element in $a^j bZ^*(G)$ for distinct
$i,j\in\{0,1,2\}$. It follows that  $|Z^*(G)|\leq 2$, otherwise
$\mathcal{E}_G$ contains a subgraph isomorphic to $K_{3,3}$.  Thus
$G$ is a non-nilpotent group of order $6$ or a non-nilpotent group
of order 12 with $|Z^*(G)|=2$. Since
$Z(A_4)=1$, the  proof of ``only if'' part is complete.\\
Conversely,  in Figure 1 we have the Engel graph of $S_3$;  the
Engel graphs of $D_{12}=\left<s,r\;|\; s^6=r^2=s^rs=1\right>$ and
$T$ are isomorphic to the following graph:
\begin{center}
\includegraphics{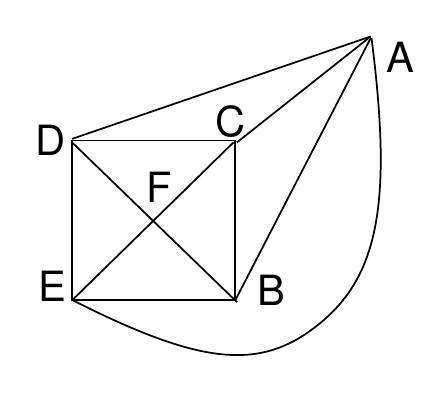}\\
Figure 2
\end{center}
where $(A,B,C,D,E,F)$ is equal to $(r,rs,rs^2,rs^4,rs^5,rs^3)$ in
$\mathcal{E}_{D_{12}}$ and $(y,yx,yx^2,yx^4,yx^5,y^3)$ in
$\mathcal{E}_T$, respectively. These graphs are visibly planar.
This completes the proof.
\end{proof}
We end this section with the following question.
\begin{qu}
Is there an infinite  non-Engel group whose Engel graph is planar?
\end{qu}
 \section{\bf Groups with the same Engel graph}
In this section we study groups with isomorphic Engel graphs. In
fact we consider the following question.
\begin{qu}\label{mainqu}  Let $G$ and $H$ be two non-Engel groups such that
$\mathcal{E}_G\cong \mathcal{E}_H$.
 For which group property $\mathcal{P}$, if  $G$ has
$\mathcal{P}$, then  $H$ also has  $\mathcal{P}$?
\end{qu}
At the moment we give the positive answer to Question
\ref{mainqu}, when $\mathcal{P}$ is the property of being finite.
\begin{thm}\label{thmfinte}
Let $G$ and $H$ be two non-Engel groups such that
$\mathcal{E}_G\cong \mathcal{E}_H$. If $G$ is finite,  then $H$ is
also a finite  group. Moreover $|L(H)|$ divides $|G|-|L(G)|$.
\end{thm}
\begin{proof}
Since $\mathcal{E}_G \cong \mathcal{E}_H$, $|H\backslash
L(H)|=|G\backslash L(G)|$. Then $H\backslash L(H)$ is finite. If
$h\in H\backslash L(H)$, then $\{h^x \;|\; x\in H\}\subseteq
H\backslash L(H)$, since $L(H)$ is closed under conjugation. Thus
every element in $H\backslash L(H)$ has finitely many conjugates
in $H$. It follows that $K=C_H(H\backslash L(H))$ has finite
index in $H$. By Corollary \ref{min-max-sol}, $\mathcal{E}_H$ has
no isolated vertex. Thus there exist two adjacent vertices $h_1$
and $h_2$ in $\mathcal{E}_H$. Now if $s\in K$, then $s\in
C_H(h_1,h_2)$. It follows that $$[sh_1,_k h_2]=[h_1,_k h_2]\not=1
\;\;\text{and}\;\; [h_2,_k sh_1]=[h_2,_k h_1]\not=1 \;\text{for
all} \; k\in\mathbb{N}.$$ Therefore $Kh_1 \subseteq H\backslash
L(H)$ and so $K$ is finite. Hence $H$ is finite and we have that
$|H|-|L(H)|=|G|-|L(G)|$. Now since $H$ is finite, it follows from
\cite{Baer} that $L(H)$ is a subgroup of $H$ and so $|L(H)|$
divides $|H|$. This completes the proof.
\end{proof}
\noindent {\bf Acknowledgements.} The author is indebted to the
referee for  his valuable comments.
 The author thanks the Centre of Excellence for Mathematics, University of Isfahan.

\end{document}